\documentclass[a4paper]{article}
\usepackage{graphics,graphicx}
\usepackage{amsfonts,amssymb,amsmath,color,mathrsfs,amstext,bm}
\usepackage{amsthm,authblk}
\usepackage{geometry}
\usepackage[colorlinks,linkcolor=red,anchorcolor=red,citecolor=red]{hyperref}
\usepackage[nameinlink,noabbrev]{cleveref}

\numberwithin{equation}{section}

\def\E{\mathbb{E}}
\newcommand{\bfP}{\mathbf{P}}
\newcommand{\R}{\mathbb{R}}

\newcommand{\dd}{\mathrm{d}}

\newcommand{\ip}[2]{\langle #1,#2\rangle}

\newtheorem{theorem}{Theorem}[section]
\newtheorem{lemma}{Lemma}[section]
\newtheorem{proposition}{Proposition}[section]

\title{On the propagation of chaos for three-dimensional Coulomb Brownian particles}

\author[a]{Lei Li\thanks{E-mail: leili2010@sjtu.edu.cn}}
\author[b]{Jian-Guo Liu\thanks{E-mail: jliu@math.duke.edu}}
\affil[a]{School of Mathematical Sciences, Institute of Natural Sciences, MOE-LSC, Shanghai Jiao Tong University, Shanghai, 200240, P.R.China.}
\affil[b]{Department of Mathematics, Department of Physics, Duke University, Durham, NC 27708, USA}
\date{}

\begin{document}
\maketitle
\begin{abstract}
While the mean field limit of the Brownian particles with 3D Coulomb interaction has 
been resolved starting from the Liouville equation, the justification 
directly from the random trajectories  is not fully resolved. In our previous work (J. Math. Phys., 60, 111501, 2019),
it has been shown that the time marginals of any limit point of the empirical measure form a weak solution 
to the mean field nonlinear Fokker-Planck equation. However, the identification
to the strong solution is not established due to the lack of regularity. In this short note,
we first make use of the variational formulation of the quadratic function to establish the 
desired space-time $L^2$ bound mentioned already in J. Math. Phys., 60, 111501, 2019.
This fills gap and shows that the time marginals must be the must be the strong solution. Moreover, making use of the virial type estimates and the force dissipation, we also establish the propagation of chaos in path space
directly from the random trajectories of the particle systems.
\end{abstract}

\section{Introduction}

Interacting particle systems arise in many problems in physics and
applied mathematics. When the number of particles is large, the mean
field approximation replaces the interaction with all other particles
by a force determined by their distribution. A basic question is
whether initially independent particles become asymptotically
independent as their number tends to infinity. This property, known
as propagation of chaos, justifies the description of a many-particle
system by a nonlinear equation for a single particle. We refer to
Kac \cite{kac1956foundations}, McKean \cite{mckean1967propagation}
and Sznitman \cite{sznitman1991topics} for the classical theory, and to
\cite{golse2016dynamics,jabin2017mean} for reviews of mean field limits.

In this note, we consider Brownian particles with repulsive Coulomb
interaction in three-dimensional space:
\begin{equation}\label{eq:particles}
 \dd X_i^N=b_i^N\,\dd t+\sqrt2\,\dd B_i,\qquad
 b_i^N=\frac1N\sum_{j\ne i}F(X_i^N-X_j^N),\qquad
 \mu_t^N=\frac1N\sum_i\delta_{X_i^N(t)}.
\end{equation}
Here $X_i^N\in\R^3$, $1\le i\le N$, the $B_i$ are independent
standard Brownian motions, and the initial positions are independent
with common probability density $\rho_0$ and independent of the
Brownian motions. The potential and force are
\begin{equation}\label{eq:kernel}
 g(x)=\frac1{4\pi|x|},\qquad
 F(x)=-\nabla g(x)=\frac{x}{4\pi|x|^3},
\end{equation}
where $-\Delta g=\delta_0$ in distributions. The factor $1/N$ gives
the mean field scaling, while the diffusion strength is fixed.

The expected limiting trajectories are independent solutions of
the McKean SDE
\begin{equation}\label{eq:mckean}
 \dd\bar X_i(t)=(F*\rho_t)(\bar X_i(t))\,\dd t+\sqrt2\,\dd B_i(t),
\end{equation}
where $\rho_t$ is the law of $\bar X_i(t)$. This density satisfies
the nonlinear Fokker--Planck equation
\begin{equation}\label{eq:nonlinearfp}
 \partial_t\rho=\Delta\rho-\nabla\cdot(b\rho)
 =\Delta\rho+\nabla\cdot(\rho\nabla h),\qquad \rho|_{t=0}=\rho_0,
\end{equation}
with
\begin{equation}\label{eq:field}
 h=g*\rho,\qquad b=F*\rho=-\nabla h.
\end{equation}
We write $\mu^N=N^{-1}\sum_i\delta_{X_i^N(\cdot)}$ for the empirical
measure on $C([0,T];\R^3)$; its time marginals are $\mu_t^N$.

There are two related approaches to the mean field limit:
one follows the trajectories and their empirical measure, while the
other studies the joint density through the Liouville equation.
For deterministic systems with Lipschitz interaction forces,
Dobrushin \cite{dobrushin1979vlasov} compared characteristic
trajectories to obtain stability in Wasserstein distance and justify
the Vlasov mean field limit.
For McKean--Vlasov diffusions with Lipschitz interactions, the
classical coupling argument compares the interacting trajectories
with independent nonlinear processes driven by the same Brownian
motions; the law of large numbers and Gr\"onwall's inequality then
give propagation of chaos
\cite{mckean1967propagation,sznitman1991topics}.
These direct trajectory estimates rely on regularity of the
interaction and require additional control for singular forces.
Along the Liouville equation approach, Jabin and Wang \cite{jabin2018quantitative} developed a relative
entropy method giving quantitative estimates with $W^{-1,\infty}$ interaction kernels. 
Bresch, Jabin and Wang
\cite{bresch2020modulated} combined relative entropy and modulated
potential energy to obtain quantitative convergence for repulsive
Coulomb and Riesz interactions for periodic systems. Their subsequent work
\cite{bresch2020mean} treats singular attractive interactions,
including the Patlak--Keller--Segel model in subcritical regimes.
On the whole space, Feng and Wang
\cite{feng2026whole,feng2026circulations} obtained quantitative
relative entropy estimates for the two-dimensional viscous vortex
model, including general circulations, while Cai, Feng, Gong and
Wang \cite{cai2024loggas} adapted the modulated free energy method
to the two-dimensional repulsive log gas.
These whole-space results give quantitative convergence of
time-marginal laws for two-dimensional systems with force
singularities. The distinction between time-marginal and path-space convergence is
essential: the former alone does not determine temporal
correlations, as emphasized in \cite{bresch2020mean}.

Regarding the three-dimensional Coulomb force of order $|x|^{-2}$ considered here,   the singularity of the force is the main difficulty in justifying this limit. For deterministic gradient
flows, Duerinckx \cite{duerinckx2016mean} treated Riesz interactions,
and Serfaty \cite{serfaty2020mean} established mean field convergence
for Coulomb and more singular Riesz interactions by modulated energy. For stochastic systems, Fournier, Hauray and Mischler
\cite{fournier2014propagation} used entropy dissipation and Fisher
information to prove path-space propagation of chaos for the
two-dimensional viscous vortex model. Liu and Yang
\cite{liu2016propagation} established noncollision for Coulomb
Brownian particles and propagation of chaos in dimension two.
The work of Bresch, Jabin and Wang
\cite{bresch2020modulated} established the propagation of chaos of Coulomb system on torus using Liouville equation for time marginals. The propagation of chaos for trajectories in full space has not been resolved as the classical coupling approach is challenging.

Our previous work \cite{LLY19} uses compactness of the empirical
path measures. Entropy and Fisher information estimates, together
with symmetrization of the interaction, show that the time marginals
of every limit point solve \eqref{eq:nonlinearfp} weakly. A
weak--strong uniqueness result was also established under an
additional space-time $L^2$ assumption on the limiting density.
However, this assumption was not justified by the available
estimates. Moreover, uniqueness of the time marginals does not
alone imply uniqueness of the probability law on path space.
Both issues must be addressed to obtain propagation of chaos
directly from the trajectories.

We resolve these two issues here. First, we retain the squared
total force in the particle energy dissipation and use its quadratic
variational formulation to pass to the limit. Combining the resulting
weighted force bound with Fisher information gives the missing
space-time $L^2$ estimate. The weak--strong uniqueness result in
\cite{LLY19} then identifies the limiting marginals with the strong
solution. Second, a localized logarithmic virial identity controls
close-particle interactions before excessive local concentration
occurs. Convergence to the bounded strong density makes this stopping
event unlikely. Together with force dissipation, this yields
convergence of the drift along the particles. Synchronous coupling
then gives convergence of each tagged trajectory in expected uniform
distance and propagation of chaos on path space, without regularizing
the particle dynamics.
The main result is the following.
\begin{theorem}\label{thm:chaos}
Assume that $\rho_0 \in H^m(\R^3)$, $m>3/2$ and has finite initial entropy, second moment and Coulomb energy
\[
\int \rho_0 \log \rho_0 \,dx<\infty, \quad \int |x|^2\rho_0\,dx<\infty, \quad \frac{1}{2}\iint g(x-y)\rho_0(x)\rho_0(y)\,dxdy<\infty.
\]
For any $T>0$, the empirical path measure $\mu^N$ converges in probability,
for the weak topology on $\bfP(C([0,T];\R^3))$, to the deterministic
law $\mu$ of the McKean SDE \eqref{eq:mckean}. Its time marginals have
the unique strong density $\rho$ solving \eqref{eq:nonlinearfp} with
initial value $\rho_0$. More precisely, use a common iid initial sequence
$X_i^N(0)=\bar X_i(0)$ and drive $X_i^N$ and $\bar X_i$ by the same
Brownian motion $B_i$. Then, for every fixed $i$,
\begin{equation}\label{eq:trajectory-convergence}
 \lim_{N\to\infty}\E\sup_{t\le T}|X_i^N(t)-\bar X_i(t)|=0.
\end{equation}
In particular, for each fixed $k$, the joint law of the first $k$
particle paths converges weakly to $\mu^{\otimes k}$ on the product
path space $C([0,T];\R^3)^k$.
\end{theorem}

The rest is organized as follows. Section \ref{eq:weakstrong} establishes the additional regularity and identifies the
limiting time marginals. Section \ref{sec:path} proves the trajectory convergence.

\section{The square integrability and weak-strong uniqueness for time marginals}\label{eq:weakstrong}

Introduce the Coulomb energy
\begin{equation}\label{eq:functionals}
 \mathcal E(\rho)=\frac12\iint g(x-y)\rho(x)\rho(y)\,\dd x\dd y,
\end{equation}
and the Fisher information
\begin{equation}\label{eq:fisher}
I(\rho)=\int \frac{|\nabla\rho|^2}{\rho}\,dx=4\int|\nabla\sqrt\rho|^2\,\dd x.
\end{equation}
For a fixed $T<\infty$, consider any subsequential limit in law of the
empirical path measures, as in \cite[Theorem 3.1]{LLY19}. 
It has been shown in \cite{LLY19} that for any subsequential limit the time marginals have densities $\mu_t(dx)=\rho_t(x)\,\dd x$. 
Note that the limit is a random. The results of \cite{LLY19} give, almost surely,
\begin{equation}\label{eq:old}
 \partial_t\rho=\Delta\rho+\nabla\cdot(\rho\nabla h),\qquad
 \int_0^T I(\rho_t)\,\dd t<\infty,\qquad
 \sup_{t\le T}\mathcal E(\rho_t)\le\mathcal E(\rho_0).
\end{equation}
They also give $\rho b\in L^1((0,T)\times\R^3)$ almost surely.
We retain the initial entropy, second moment, and finite-energy
hypotheses of that theorem throughout.

\subsection{Mean field energy dissipation}

The key result we establish in this note will be the following mean field energy dissipation
anticipated in \cite[Remark 4.2]{LLY19}.
\begin{theorem}\label{thm:main}
Under the settings above, every empirical-measure limit satisfies
\begin{equation}\label{eq:new}
 \mathbb E\int_0^T\int\rho|\nabla h|^2\,\dd x\dd t
 \le\mathcal E(\rho_0).
\end{equation}
\end{theorem}

 We prove it directly from the energy dissipation involving forces from the particle estimate. Note that this cannot be obtained by testing the limiting equation with $h$, which is not available by \eqref{eq:old}.
 The key point is to make use of the energy dissipation of the particle system, formulated in a variational form.
 In \cite{LLY19}, this energy dissipation was not fully made use of due to the nonlinearity.

\begin{proof}[Proof of Theorem \ref{thm:main}]
Define the energy associated with the particle system by
\begin{equation}
 \mathcal E_N(X)=\frac1{2N^2}\sum_{i\ne j}g(X_i-X_j).
\end{equation}
It has been justified in \cite[Equation (2.49)]{LLY19} that
\begin{equation}\label{eq:DN}
 \mathbb E D_N\le\mathbb E\mathcal E_N(X(0))
 =\frac{N-1}{N}\mathcal E(\rho_0),\qquad
 D_N:=\int_0^T\frac1N\sum_i|b_i^N|^2\,\dd t.
\end{equation}

To make use of the dissipation $D_N$, the key is to use the quadratic duality
\[
|b|^2=\sup_{a}(2a\cdot b-|a|^2).
\]
Fix any particular limit $\rho_t$. Let $m(\dd t\dd x)=\rho_t(x)\,\dd t\dd x$.
It is a finite measure, and $b\in L^1(m)$. Quadratic duality gives
\begin{equation}\label{eq:duality}
 \int|b|^2\,\dd m
 =\sup_{a\in C_c^\infty}\left\{2\int a\cdot b\,\dd m
                           -\int|a|^2\,\dd m\right\}=:\sup_{a\in C_c^{\infty}}(2J(a)-Q(a)),
\end{equation}
including when the left side is infinite. Note that symmetrization gives
\begin{equation}\label{eq:Jlimit}
 J(a)=\int_0^T\int a\cdot b\rho\,\dd x\dd t
 =\frac12\int_0^T\iint
 [a(t,x)-a(t,y)]\cdot F(x-y)\rho_t(x)\rho_t(y)
 \,\dd x\dd y\dd t.
\end{equation}
Correspondingly, we can introduce for smooth vector field $a\in C_c^\infty((0,T)\times\R^3;\R^3)$ that
\begin{equation}
 J_N(a)=\int_0^T\frac1N\sum_i a(t,X_i^N)\cdot b_i^N\,\dd t,\quad
 Q_N(a)=\int_0^T\ip{\mu_t^N}{|a(t)|^2}\,\dd t.
\end{equation}
The oddness of $F$ gives again the symmetrization
\begin{equation}\label{eq:sym}
 J_N(a)=\frac1{2N^2}\int_0^T\sum_{i\ne j}
 [a(t,X_i^N)-a(t,X_j^N)]\cdot F(X_i^N-X_j^N)\,\dd t.
\end{equation}
The singularity in the integrand of $J$ is
only $|x-y|^{-1}$, exactly as in the weak-equation argument of
\cite{LLY19}.

Choose a countable family $\{a_k\}_{k\ge1}$ of smooth compactly
supported vector fields including zero function, dense in the uniform norm
on any compact set.  For every finite $K$, one clearly has
\begin{equation}\label{eq:finite}
 \max_{1\le k\le K}\{2J_N(a_k)-Q_N(a_k)\}\le D_N.
\end{equation}
To pass the limit, we choose a cutoff parameter $\delta$ to define $g_{\delta}, F_{\delta}$ and thus $J^{\delta}, J_N^{\delta}$
following the same notation as in \cite[Equation (2.24)]{LLY19}.
For a fixed cutoff $\delta$, the left side is bounded and continuous as a
functional of the empirical path measure, so its expectation
converges as $N\to\infty$:
\[
\lim_{N\to\infty}\E  \max_{1\le k\le K}\{2J_N^{\delta}(a_k)-Q_N(a_k)\}= \E \max_{1\le k\le K}\{2J^{\delta}(a_k)-Q(a_k)\}.
\]
Later, in Lemma \ref{lem:cutoff}, we will show that for any fixed $a$ and $\alpha\in (1, 2)$.
\begin{equation}
 \sup_N\mathbb E|J_N(a)-J_N^\delta(a)|
 +\mathbb E|J(a)-J^\delta(a)|
 \le C_{a,T,\alpha}\delta^{\alpha-1},
\end{equation}
Combing these results, since there are only finitely many $a_k$ for fixed $K$, taking first $N\to\infty$ and then $\delta\to 0$ yields
\begin{equation}\label{eq:finite-limit}
 \mathbb E\max_{k\le K}\{2J(a_k)-Q(a_k)\}  \le\liminf_N\mathbb E D_N \le \mathcal E(\rho_0).
\end{equation}
The liminf in $\E D_N$ is realized if  one restricts to a further subsequence attaining it.

Because zero is included, the maxima in \eqref{eq:finite-limit}
are nonnegative and increase with $K$. Monotone convergence then proves 
\[
\lim_{K\to\infty}\mathbb E\max_{k\le K}\{2J(a_k)-Q(a_k)\}
=\mathbb E\max_{k}\{2J(a_k)-Q(a_k)\}=\mathbb E\int_0^T\int\rho|b|^2\,\dd x\dd t.
\]
The last equality holds because $\{a_k\}$ is dense for continuous functions on any compact sets. 
One may restrict $b$ on compact sets and use the uniform approximation of
$b\in L^1(m)$ using continuous functions. Passing the limit on the truncation gives the claim.
Hence, one eventually has
\begin{equation}
 \mathbb E\int_0^T\int\rho|b|^2\,\dd x\dd t
 \le\liminf_{N\to\infty}\mathbb E D_N\le \mathcal E(\rho_0).
\end{equation}
\end{proof}

\begin{lemma}\label{lem:cutoff}
For every fixed smooth vector test $a$, there is some $\alpha\in(1,2)$ such that
\begin{equation}\label{eq:cutoff}
 \sup_N\mathbb E|J_N(a)-J_N^\delta(a)|
 +\mathbb E|J(a)-J^\delta(a)|
 \le C_{a,T,\alpha}\delta^{\alpha-1},
\end{equation}
\end{lemma}

\begin{proof}
Since $F_{\delta}$ is continuous and bounded, $J_N^{\delta}$ and $J^{\delta}$ are  bounded continuous
functional for the weak topology used in \cite{LLY19}. The same holds for $Q_N$.
Now, we find that
\[
|J_N(a)-J_N^{\delta}(a)|\le \frac{1}{2N^2}\int_0^T \sum_{i\neq j}
|a(t,X_i^N)-a(t,X_j^N)| |F(X_i^N-X_j^N)-F_{\delta}(X_i^N-X_j^N)|\,\dd t.
\]
Since $a(t,x)-a(t,y)$ gives one $|x-y|$, one finds that
\[
|a(t,X_i^N)-a(t,X_j^N)| |F(X_i^N-X_j^N)-F_{\delta}(X_i^N-X_j^N)|\le
C_a \delta^{\alpha-1} |X_i^N-X_j^N|^{-\alpha} 1_{|X_i^N-X_j^N| \le \delta}.
\]

Let $f_t^{N,(2)}$ denote the two-particle marginal density. It has been shown in \cite{LLY19} that
\begin{equation}\label{eq:pairfisher}
 \sup_{N\ge2}\int_0^T I(f_t^{N,(2)})\,\dd t<\infty.
\end{equation}
The inverse-distance estimate of \cite[Lemma 2.5]{LLY19} gives
\begin{equation}\label{eq:inverse}
 \sup_{N\ge2}\mathbb E\int_0^T\frac1{N^2}
 \sum_{i\ne j}|X_i^N-X_j^N|^{-\alpha}\,\dd t\le C_\alpha(1+I(f_t^{N,(2)})),
 \qquad 1<\alpha<2.
\end{equation}
Then, one obtains the particle part of \eqref{eq:cutoff}.
For the limit, the same estimate follows by applying the
inverse-distance bound to $\rho_t\otimes\rho_t$, whose Fisher
information is $2I(\rho_t)$, and using
$\mathbb E\int_0^T I(\rho_t)\,\dd t<\infty$, already proved in
\cite{LLY19}. The symmetrized integral is therefore absolutely
integrable. The unsymmetrized expression in \eqref{eq:Jlimit} is
also legitimate by $\rho b\in L^1$ and the absolute convolution
integrability established there; oddness then identifies the two.
\end{proof}

\subsection{The square integrability of $\rho$}

Here, we establish the desired control in \cite{LLY19} and conclude the weak-strong uniqueness.
\begin{proposition}\label{pro:weakstrong}
For the density of the time marginals of any weak limit, it holds that
\[
\E  \int_0^T\int\rho^2\,\dd x\dd t<\infty.
\]
Consequently, $\rho_t$ must be the unique strong solution to the nonlinear Fokker-Planck equation 
 \eqref{eq:nonlinearfp}.
\end{proposition}

\begin{lemma}\label{lem:spatial}
Suppose $\rho\ge0$ has mass one for almost every time,
$\int_0^T I(\rho_t)\,\dd t<\infty$, and $h=g*\rho$ satisfies
$\int_0^T\int\rho|\nabla h|^2<\infty$. Then
\begin{equation}\label{eq:L2}
 \int_0^T\int\rho^2\,\dd x\dd t
 \le\left(\int_0^T I(\rho_t)\,\dd t\right)^{1/2}
       \left(\int_0^T\int\rho|\nabla h|^2\,\dd x\dd t\right)^{1/2}.
\end{equation}
\end{lemma}

\begin{proof}
For almost every fixed $t$, $\sqrt\rho\in H^1(\R^3)$, so that
$\rho\in L^1\cap L^3$ and
$\nabla\rho=2\sqrt\rho\nabla\sqrt\rho\in L^{3/2}$.
In particular, $\rho\in L^2$ at that time.

Elliptic estimates give $\nabla h\in W^{1,3}_{\mathrm{loc}}$ and
$\Delta h=-\rho$. The local Sobolev product rule is legitimate
(for example use $\nabla h\in L^3_{\mathrm{loc}}$), and gives
\begin{equation}\label{eq:product}
 \nabla\cdot(\rho\nabla h)=\nabla\rho\cdot\nabla h-\rho^2
 \quad\hbox{in spatial distributions for almost every }t.
\end{equation}
Choose $\chi_R(x)=\chi(x/R)$ with $0\le\chi_R\uparrow1$,
$\chi\in C_c^\infty$, and $|\nabla\chi_R|\le C/R$.
Testing \eqref{eq:product} gives
\begin{equation}\label{eq:cutidentity}
 \int\chi_R\rho^2
 =\int\chi_R\nabla\rho\cdot\nabla h
   +\int\rho\nabla h\cdot\nabla\chi_R.
\end{equation}
The first term is absolutely integrable in space-time because
\begin{equation}\label{eq:cross}
 \int_0^T\int|\nabla\rho\cdot\nabla h|
 \le2\left(\int_0^T\int|\nabla\sqrt\rho|^2\right)^{1/2}
       \left(\int_0^T\int\rho|\nabla h|^2\right)^{1/2}.
\end{equation}
The boundary term integrated in time (the last term in in \eqref{eq:cutidentity}) is at most
\begin{equation}
 \frac C R\int_0^T\int\rho|\nabla h|
 \le\frac{C\sqrt T}R
       \left(\int_0^T\int\rho|\nabla h|^2\right)^{1/2}.
\end{equation}
Integrate \eqref{eq:cutidentity} in time and send $R\to\infty$.
Monotone convergence on its nonnegative left side and dominated
convergence on its first right side shows that
\begin{equation}\label{eq:exact}
 \int_0^T\int\rho^2
 =\int_0^T\int\nabla\rho\cdot\nabla h.
\end{equation}
Equation \eqref{eq:cross} then gives \eqref{eq:L2}.
\end{proof}

We now conclude the main claim in this section.
\begin{proof}[Proof of Proposition \ref{pro:weakstrong}]
Writing $C_I(T)=\mathbb E\int_0^T I(\rho_t)\,\dd t$, the Cauchy--Schwarz on the probability space yield
\begin{equation}\label{eq:expectedL2}
 \mathbb E\int_0^T\int\rho^2
 \le\sqrt{C_I(T)\mathcal E(\rho_0)}<\infty.
\end{equation}
It supplies the additional estimate appearing in the
statement of \cite[Theorem 4.1]{LLY19}.

By \eqref{eq:expectedL2}, every subsequential limit of the empirical
path measures has marginal densities in $L^2(0,T;L^2)$ almost surely. Applying
\cite[Proposition 4.1]{LLY19} identifies all these marginals with the
same deterministic strong solution $\rho_t$. 
\end{proof}

\section{The propagation of chaos in path space}\label{sec:path}

Identifying the
marginals alone does not identify a path law. We complete the last step for the propagation of chaos in path space,
using a stopped logarithmic virial estimate and synchronous coupling.

\begin{proof}[Proof of Theorem~\ref{thm:chaos}]

Part of the assertions have been proved in Proposition \ref{pro:weakstrong}. Here, we focus on the path convergence.

Write $d_{\mathrm{BL}}$ for the bounded Lipschitz distance on
$\bfP(\R^3)$. Tightness of the empirical path measures, already
proved in \cite{LLY19}, and marginal identification imply
\begin{equation}\label{eq:spatial-curve-convergence}
 \sup_{t\le T}d_{\mathrm{BL}}(\mu_t^N,\rho_t\,dx)
 \longrightarrow0\quad\hbox{in probability}.
\end{equation}
Indeed, the map from a path law to its curve of time marginals is
continuous from the weak topology on $\bfP(C([0,T];\R^3))$ to the
uniform $d_{\mathrm{BL}}$ topology. Note that we have not identified the limiting path law yet at this point.
The strong regularity $\rho\in C([0,T];H^m)$, $m>3/2$, gives
\begin{equation}\label{eq:strong-field-bounds}
 M:=\sup_{t\le T}\|\rho_t\|_\infty<\infty,\qquad
 B:=\sup_{t\le T}\|b_t\|_\infty<\infty,\qquad
 L_T:=\int_0^T\|\nabla b_t\|_\infty\,dt<\infty.
\end{equation}
Here $b_t=F*\rho_t$. The bound on $b$ follows by splitting the
convolution at unit distance. For its derivative, $\nabla b$ is a
matrix of order-zero Fourier multipliers applied to $\rho$, bounded
in $H^m$; Sobolev embedding applies.

Next, we use the synchronous coupling.  Construct $\bar X_i$ with the same initial position and Brownian
motion as $X_i^N$, and deterministic drift $b_t=F*\rho_t$. The
processes $\bar X_i$ are i.i.d.. Their time marginals are $\rho_t$ by
uniqueness for the linear Fokker--Planck equation with this bounded
Lipschitz drift, so they solve \eqref{eq:mckean}.
Define
\[
 \mathcal R_N=\int_0^T\frac1N\sum_i
       |b_i^N(t)-b_t(X_i^N(t))|\,dt.
\]
Subtracting the SDEs and using \eqref{eq:strong-field-bounds},
Gr\"onwall's inequality gives
\[
 \sup_{t\le T}|X_i^N(t)-\bar X_i(t)|
 \le e^{L_T}\int_0^T
             |b_i^N(t)-b_t(X_i^N(t))|\,dt.
\]
By exchangeability, it thus holds that
\[
\E \sup_{t\le T}|X_i^N(t)-\bar X_i(t)| \le C(T)\E\mathcal R_N.
\]
Choose the cutoff force $F^\delta$ again so that
\[
|F-F^\delta|\le C|z|^{-2}\mathbf1_{\{|z|<\delta\}}.
\]
For any stopping time $0\le\tau\le T$, define the random variable
\begin{gather}\label{eq:Udelta}
 U_N(\delta;\tau):=\int_0^\tau\frac1{N^2}\sum_{i\ne j}
 \frac{\mathbf1_{\{|X_i^N(t)-X_j^N(t)|<\delta\}}}
 {|X_i^N(t)-X_j^N(t)|^2}\,dt,\qquad
 U_N(\delta):=U_N(\delta;T).
\end{gather}
Thus $U_N(\delta)$ is the full-time random integral; the expectation
of the stopped integral will be written $\E U_N(\delta;\theta_N^\delta)$.
Splitting the drift error yields
\begin{align}
 \mathcal R_N\le{}& C U_N(\delta)
 +T\sup_{t\le T,x\in\R^3}
       |\frac{1}{N}\sum_{j\neq i}F^\delta(x-X_j^N)-(F^\delta*\rho_t)(x)|
 +CTM\delta.
 \label{eq:drift-error-decomposition}
\end{align}
For fixed $\delta$, the middle term tends to zero in probability
by \eqref{eq:spatial-curve-convergence}, since translates of
$F^\delta$ have a common bounded Lipschitz norm. The last term uses
$\int_{|z|<\delta}|F(z)|\,dz\le C\delta$.
The main technical difficulty in this part is to control $U_N(\delta)$ in probability, 
which then implies that $\mathcal R_N\to 0$ in probability as $N\to\infty$ so that $\E \mathcal R_N\to 0$.

If we apply the three-dimensional Hardy inequality, applied in the first variable of the two-particle density, one may use the Fisher information of the two-particle distribution to control the following
\[
 \iint\frac{f_t^{N,(2)}(x,y)}{|x-y|^2}\,dxdy
 \le4\iint|\nabla_x\sqrt{f_t^{N,(2)}}|^2\,dxdy.
\]
Consequently, \eqref{eq:pairfisher} and exchangeability yield
\begin{equation}\label{eq:hardy-endpoint}
 C_H:=\sup_{N\ge2}\E\int_0^T\frac1{N^2}
 \sum_{i\ne j}|X_i^N-X_j^N|^{-2}\,dt<\infty.
\end{equation}
This bound alone does not give uniform integrability near
the diagonal and it does not give $\E U_N(\delta)\to 0$ as $\delta\to 0$.
 We will use this for the control before a local concentration stopping time. To get the desired estimate, we will use 
 a virial type theorem.

 The observation is that the quadratic variation term of $\log |x|$ will yield
 $1/|x|^2$, while the virial type theorem will yield symmetrization on the force term so that it may be controlled.

 Choose a smooth radial cutoff $0\le\chi\le1$, equal to one on $B_1$ and zero outside $B_2$. For $z\ne0$, set
\[
 \psi_\delta(z)=\chi(z/\delta)\log\frac{|z|}{2\delta},
 \qquad a_\delta(z)=\nabla\psi_\delta(z).
\]
Since $\Delta\log|z|=|z|^{-2}$ in three dimensions,
\begin{gather}\label{eq:log-cutoff-properties}
 |\psi_\delta(z)|\le C\delta g(z),\qquad
 |a_\delta(z)|\le\frac C{|z|}\mathbf1_{\{|z|<2\delta\}},\\
 \Delta\psi_\delta(z)
 =\frac{\mathbf1_{\{|z|<\delta\}}}{|z|^2}+R_\delta(z),
 \qquad
 |R_\delta(z)|\le\frac C{\delta^2}
       \mathbf1_{\{\delta\le|z|\le2\delta\}}.\nonumber
\end{gather}
Introduce
\[
 W_N^\delta=\frac1{2N^2}\sum_{i\ne j}
                 \psi_\delta(X_i^N-X_j^N),\qquad
 A_i^\delta=\frac1N\sum_{j\ne i}a_\delta(X_i^N-X_j^N).
\]
To treat the singularity at $z=0$, for fixed $N$ and $\delta$, we set
\[
 \sigma_\varepsilon=\inf\left\{t\in[0,T]:
 \min_{i<j}|X_i^N(t)-X_j^N(t)|\le\varepsilon\right\}\wedge T,
 \qquad \inf\varnothing=+\infty.
\]
The finite-$N$ system is almost surely noncolliding. Since its
trajectories are continuous and there are finitely many pairs,
\[
 \min_{t\in[0,T],\,i<j}|X_i^N(t)-X_j^N(t)|>0
 \quad\hbox{almost surely}.
\]
Consequently, $\sigma_\varepsilon=T$ for all sufficiently small
$\varepsilon$, almost surely.

The function $W_N^\delta$ is $C^2$ on the open
collision-free configuration space. The localized It\^o formula
therefore applies up to $\sigma_\varepsilon$. Using
$\nabla_{X_i}W_N^\delta=A_i^\delta/N$ and the oddness of $a_\delta$
gives, initially up to this stopping time,
\begin{equation}\label{eq:log-virial}
 dW_N^\delta=
 \left[\frac1N\sum_i A_i^\delta\cdot b_i^N
 +\frac1{N^2}\sum_{i\ne j}\Delta\psi_\delta(X_i^N-X_j^N)\right]dt
 +d\mathcal M_N^\delta.
\end{equation}
Here $d\mathcal M_N^\delta=(\sqrt2/N)\sum_i A_i^\delta\cdot dB_i$.

Jensen's inequality
and \eqref{eq:hardy-endpoint} imply that one can take $\varepsilon\to 0$ in expectation and obtain
\[
 \E\int_0^T\frac1N\sum_i|A_i^\delta|^2\,dt\le C C_H,
 \qquad
 \E\langle\mathcal M_N^\delta\rangle_T
 =\frac2N\E\int_0^T\frac1N\sum_i|A_i^\delta|^2\,dt<\infty.
\]
Thus $\mathcal M_N^\delta$ is a square-integrable martingale. The two finite-variation terms are absolutely integrable, since
\begin{align*}
 \E\int_0^T\frac1N\sum_i|A_i^\delta\cdot b_i^N|\,dt
 &\le (C C_H\E D_N)^{1/2}<\infty,\\
 \E\int_0^T\frac1{N^2}\sum_{i\ne j}
 |\Delta\psi_\delta(X_i^N-X_j^N)|\,dt
 &\le C C_H+CT\delta^{-2}<\infty.
\end{align*}
Hence,  for any bounded stopping time $\theta\le T$, It\^o isometry and
dominated convergence give
\[
 \mathcal M_N^\delta(\theta\wedge\sigma_\varepsilon)
 \longrightarrow\mathcal M_N^\delta(\theta)
 \quad\hbox{in }L^2.
\]
The expectation of the two $dt$ terms also have $L^1$ convergence when the limit $\varepsilon\to 0$ is taken on the endpoint
$\theta \wedge \sigma_{\varepsilon}$. Then, one has the identity
\begin{equation}\label{eq:stoppedvirial}
 \E[W_N^\delta(\theta)-W_N^\delta(0)]
 =  \E\int_0^\theta\frac1N\sum_i A_i^\delta\cdot b_i^N\,dt +\E\int_0^\theta\frac1{N^2}\sum_{i\ne j}
       \Delta\psi_\delta(X_i^N-X_j^N)\,dt.
\end{equation}

We expect that the particles will not gather together at some point. For this purpose, 
choose $w\in C_c^\infty(\R^3)$ with $0\le w\le1$, $w=1$ on $B_2$,
and $w=0$ outside $B_3$. For $0<\delta<1$, define
\[
 c_N^\delta(t)=\sup_{x\in\R^3}
       \int w((x-y)/\delta)\,\mu_t^N(dy),\qquad
 m_\delta=2M\|w\|_{L^1}\delta^3.
\]
This is the fraction of particles near $x$. We will take $\theta=\theta_N^{\delta}$, defined by
\[
 \theta_N^\delta=\inf\{t\in[0,T]:c_N^\delta(t)\ge m_\delta\}\wedge T,
 \qquad \inf\varnothing=+\infty.
\]
The equation \eqref{eq:spatial-curve-convergence} gives
uniform convergence in probability of $c_N^\delta$ to its
deterministic counterpart as $N\to\infty$
\[
\sup_x\int w(\frac{x-y}{\delta})\rho_t\,dx\le  M\|w\|_{L^1}\delta^3=m_\delta/2.
\]
 In particular,
\begin{equation}\label{eq:rare-stop}
 \lim_{N\to\infty}\bfP(\theta_N^\delta<T)=0
 \qquad\hbox{for every fixed }\delta>0.
\end{equation}
Before $\theta=\theta_N^\delta$, every ball of radius $2\delta$ has empirical
mass at most $m_\delta$.

Before $\theta$, Cauchy--Schwarz in the
particle index gives the crucial estimate
\begin{align*}
 |A_i^\delta|^2
 &\le C\left(\frac1N\sum_{j\ne i}
          \mathbf1_{\{|X_i^N-X_j^N|<2\delta\}}\right)
       \left(\frac1N\sum_{j\ne i}
          \frac{\mathbf1_{\{|X_i^N-X_j^N|<2\delta\}}}
               {|X_i^N-X_j^N|^2}\right)\\
 &\le C m_\delta\,\frac1N\sum_{j\ne i}|X_i^N-X_j^N|^{-2}.
\end{align*}
Thus, \eqref{eq:hardy-endpoint} yields that
\begin{equation}\label{eq:stopped-K}
 \E\int_0^\theta\frac1N\sum_i|A_i^\delta|^2\,dt
 \le C m_\delta C_H.
\end{equation}
The force dissipation ($ \E D_N$) now yields
\begin{equation}\label{eq:stopped-virial-drift}
 \E\left|\int_0^\theta\frac1N\sum_i A_i^\delta\cdot b_i^N\,dt\right|
 \le C\sqrt{m_\delta C_H\mathcal E(\rho_0)}.
\end{equation}

By \eqref{eq:log-cutoff-properties}, $|W_N^{\delta}|$ is controlled by $\mathcal E_N$. Using the energy dissipation in \cite[Equation (2.49)]{LLY19}, one has
\begin{equation}\label{eq:log-endpoints}
 \E|W_N^\delta(\theta)|+\E|W_N^\delta(0)|
 \le C\delta\mathcal E(\rho_0).
\end{equation}
Before $\theta$, the averaged annulus mass is at most $m_\delta$.
The integrated identity \eqref{eq:log-virial} and
\eqref{eq:log-cutoff-properties} give (recalling \eqref{eq:Udelta} for $U_N(\delta;\theta)$)
\begin{align*}
 \E U_N(\delta;\theta)
 ={}&\E[W_N^\delta(\theta)-W_N^\delta(0)]
 -\E\int_0^\theta\frac1N\sum_i A_i^\delta\cdot b_i^N\,dt\\
 &-\E\int_0^\theta\frac1{N^2}\sum_{i\ne j}
                   R_\delta(X_i^N-X_j^N)\,dt.
\end{align*}
Bounding the three terms separately by
\eqref{eq:log-endpoints}, \eqref{eq:stopped-virial-drift}, and
$CTm_\delta/\delta^2$, respectively, proves
\begin{equation}\label{eq:stopped-near-force}
 \E U_N(\delta;\theta) \le C\left[
 \delta\mathcal E(\rho_0)+\frac{Tm_\delta}{\delta^2}
 +\sqrt{m_\delta C_H\mathcal E(\rho_0)}\right] \le C_T(\delta+\delta^{3/2}).
\end{equation}
All constants here are independent of $N$ and $0<\delta<1$.

On the event $\{\theta_N^\delta=T\}$, the stopped and full integrals
coincide. Therefore, for every $\eta>0$, Markov's inequality gives
\begin{align*}
 \bfP(U_N(\delta)>\eta)
 &\le\bfP(\theta_N^\delta<T)
       +\bfP(U_N(\delta;\theta_N^\delta)>\eta)\\
 &\le\bfP(\theta_N^\delta<T)
       +\frac{C_T}{\eta}(\delta+\delta^{3/2}).
\end{align*}
Consequently,
\begin{equation}\label{eq:near-force-in-probability}
 \lim_{\delta\downarrow0}\limsup_{N\to\infty}
       \bfP(U_N(\delta)>\eta)=0.
\end{equation}
Consequently, 
\[
\limsup_{N\to\infty}
       \bfP(\mathcal{R}_N>\eta)=0.
\]

By Cauchy--Schwarz and \eqref{eq:DN},
\[
 \mathcal R_N\le\sqrt{TD_N}+TB,\qquad
 K_T:=\sup_N\E\mathcal R_N^2
 \le2T\mathcal E(\rho_0)+2T^2B^2<\infty.
\]
For every $\eta>0$, split the expectation according to whether
$\mathcal R_N$ exceeds $\eta$. Cauchy--Schwarz gives
\begin{align*}
 \E\mathcal R_N
 &=\E[\mathcal R_N\mathbf1_{\{\mathcal R_N\le\eta\}}]
   +\E[\mathcal R_N\mathbf1_{\{\mathcal R_N>\eta\}}]\\
 &\le\eta+
   (\E\mathcal R_N^2)^{1/2}\bfP(\mathcal R_N>\eta)^{1/2}.
\end{align*}
Taking $N\to\infty$
yields $\limsup_{N\to\infty}\E\mathcal R_N\le\eta$.
Since $\eta$ is arbitrary, one has
\begin{equation}\label{eq:drift-error-L1}
 \E\mathcal R_N\longrightarrow0.
\end{equation}

 Equation \eqref{eq:drift-error-L1} proves \eqref{eq:trajectory-convergence}. Summing over any fixed number of
labels gives convergence to independent McKean paths.  Finally, use the bounded uniform metric
$1\wedge\|\omega-\eta\|_\infty$ on path space. The bounded Lipschitz
distance between $\mu^N$ and
$N^{-1}\sum_i\delta_{\bar X_i(\cdot)}$ is at most
\[
 \frac1N\sum_i\left(1\wedge
             \sup_{t\le T}|X_i^N(t)-\bar X_i(t)|\right),
\]
whose expectation tends to zero. The empirical measure of the iid
$\bar X_i$ converges almost surely to their common path law $\mu$.
This proves the asserted empirical path convergence in probability.
\end{proof}

\bibliographystyle{plain}
\bibliography{meanfield}

\begin{thebibliography}{10}

\bibitem{bresch2020modulated}
Didier Bresch, Pierre-Emmanuel Jabin, and Zhenfu Wang.
\newblock Modulated free energy and mean field limit.
\newblock {\em S{\'e}minaire Laurent Schwartz—EDP et applications}, pages
  1--22, 2020.

\bibitem{bresch2020mean}
Didier Bresch, Pierre-Emmanuel Jabin, and Zhenfu Wang.
\newblock Mean field limit and quantitative estimates with singular attractive
  kernels.
\newblock {\em Duke Mathematical Journal}, 172(13):2591--2641, 2023.

\bibitem{cai2024loggas}
Shuzhe Cai, Xuanrui Feng, Yun Gong, and Zhenfu Wang.
\newblock Propagation of chaos for {2D} log gas on the whole space.
\newblock {\em arXiv preprint arXiv:2411.14777}, 2024.

\bibitem{dobrushin1979vlasov}
Roland~L'vovich Dobrushin.
\newblock Vlasov equations.
\newblock {\em Funktsional'nyi Analiz i ego Prilozheniya}, 13(2):48--58, 1979.

\bibitem{duerinckx2016mean}
Mitia Duerinckx.
\newblock Mean-field limits for some {Riesz} interaction gradient flows.
\newblock {\em SIAM Journal on Mathematical Analysis}, 48(3):2269--2300, 2016.

\bibitem{feng2026whole}
Xuanrui Feng and Zhenfu Wang.
\newblock Quantitative propagation of chaos for {2D} viscous vortex model on
  the whole space.
\newblock {\em Peking Mathematical Journal}, pages 1--31, 2026.

\bibitem{feng2026circulations}
Xuanrui Feng and Zhenfu Wang.
\newblock Quantitative propagation of chaos for {2D} viscous vortex model with
  general circulations on the whole space.
\newblock {\em Nonlinearity}, 39(1):015009, 2026.

\bibitem{fournier2014propagation}
Nicolas Fournier, Maxime Hauray, and St{\'e}phane Mischler.
\newblock Propagation of chaos for the {2D} viscous vortex model.
\newblock {\em Journal of the European Mathematical Society}, 16(7):1423--1466,
  2014.

\bibitem{golse2016dynamics}
Fran{\c{c}}ois Golse.
\newblock On the dynamics of large particle systems in the mean field limit.
\newblock In {\em Macroscopic and large scale phenomena: coarse graining, mean
  field limits and ergodicity}, pages 1--144. Springer, 2016.

\bibitem{jabin2017mean}
Pierre-Emmanuel Jabin and Zhenfu Wang.
\newblock Mean field limit for stochastic particle systems.
\newblock In {\em Active Particles, Volume 1: Advances in Theory, Models, and
  Applications}, pages 379--402. Birkh{\"a}user, 2017.

\bibitem{jabin2018quantitative}
Pierre-Emmanuel Jabin and Zhenfu Wang.
\newblock Quantitative estimates of propagation of chaos for stochastic systems
  with {$W^{-1,\infty}$} kernels.
\newblock {\em Inventiones mathematicae}, 214:523--591, 2018.

\bibitem{kac1956foundations}
Mark Kac.
\newblock Foundations of kinetic theory.
\newblock In {\em Proceedings of the Third Berkeley Symposium on Mathematical
  Statistics and Probability}, volume~3, pages 171--197. University of
  California Press, 1956.

\bibitem{LLY19}
Lei Li, Jian-Guo Liu, and Pu~Yu.
\newblock On the mean field limit for {Brownian} particles with {Coulomb}
  interaction in {3D}.
\newblock {\em Journal of Mathematical Physics}, 60(11):111501, 2019.

\bibitem{liu2016propagation}
Jian-Guo Liu and Rong Yang.
\newblock Propagation of chaos for large {Brownian} particle system with
  {Coulomb} interaction.
\newblock {\em Research in the Mathematical Sciences}, 3:40, 2016.

\bibitem{mckean1967propagation}
Henry~P. McKean, Jr.
\newblock Propagation of chaos for a class of non-linear parabolic equations.
\newblock In {\em Stochastic Differential Equations}, Lecture Series in
  Differential Equations, pages 41--57. Air Force Office of Scientific
  Research, Arlington, Virginia, 1967.
\newblock Session 7, Catholic University.

\bibitem{serfaty2020mean}
Sylvia Serfaty.
\newblock Mean field limit for {Coulomb}-type flows.
\newblock {\em Duke Mathematical Journal}, 169(15):2887--2935, 2020.
\newblock With an appendix by Mitia Duerinckx and Sylvia Serfaty.

\bibitem{sznitman1991topics}
Alain-Sol Sznitman.
\newblock Topics in propagation of chaos.
\newblock In {\em {\'E}cole d'{\'E}t{\'e} de Probabilit{\'e}s de Saint-Flour
  XIX---1989}, volume 1464 of {\em Lecture Notes in Mathematics}, pages
  165--251. Springer, 1991.

\end{thebibliography}

\end{document}